\documentclass{elsarticle}
\usepackage{amsmath}
\usepackage{amssymb}
\usepackage{amsfonts}
\usepackage{rotating}
\usepackage{graphicx}
\usepackage{floatflt,epsfig}
\usepackage{lineno,hyperref}
\usepackage{enumerate}
\usepackage{colortbl}
\usepackage{array,tabularx,tabulary,booktabs}
\usepackage{longtable}
\usepackage{multirow}
\usepackage{wrapfig}
\usepackage{subcaption}
\newcolumntype{^}{>{\currentrowstyle}}

\journal{arXiv}

\def\ms{\medskip}
\def\bs{\bigskip}
\def\st{:}
\def\eps{\varepsilon}

\newtheorem{conjecture}{Conjecture}
\newtheorem{theorem}{Theorem}
\newcommand{\proof}{\medskip\noindent{\bf Proof.~}}

\begin{document}
\renewcommand{\abstractname}{Abstract}
\renewcommand{\refname}{References}
\renewcommand{\tablename}{Table.}
\renewcommand{\arraystretch}{0.9}
\thispagestyle{empty}
\sloppy

\begin{frontmatter}
\title{An improved bound on the chromatic number of the Pancake graphs}

\author{Leen Droogendijk}
\ead{drooge001@kpnmail.nl}

\author[01,02]{Elena~V.~Konstantinova}
\ead{e\_konsta@math.nsc.ru}

\address[01]{Sobolev Institute of Mathematics, Ak. Koptyug av. 4, Novosibirsk, 630090, Russia}
\address[02]{Novosibisk State University, Pirogova str. 2, Novosibirsk, 630090, Russia}

\begin{abstract}
In this paper an improved bound on the chromatic number of the Pancake graph $P_n, n\geqslant 2$, is presented. The bound is obtained using a subadditivity property of the chromatic number of the Pancake graph. We also investigate an equitable coloring of $P_n$. An equitable $(n-1)$-coloring based on efficient dominating sets is given and optimal equitable $4$-colorings are considered for small $n$. It is conjectured that the chromatic number of $P_n$ coincides with its equitable chromatic number for any $n\geqslant 2$.
\end{abstract}

\begin{keyword} Pancake graph; chromatic number; equitable coloring

\vspace{\baselineskip}
\MSC[2010] 05C15\sep 05C25\sep 05C69
\end{keyword}
\end{frontmatter}

\section{Introduction}

The {\it Pancake graph} $P_n, \ n\geqslant 2$, is defined as the Cayley graph over the symmetric group $\mathrm{Sym}_n$ with the generating set of all prefix--reversals $r_i, \ 2\leqslant i \leqslant n$, inverting the order of any substring $[1,i]$ of a permutation when multiplied on the right. It is a connected vertex--transitive $(n-1)$-regular graph without loops and multiple edges of order $n!$. It contains all cycles $C_l$ of length $l$, where $6\leqslant l \leqslant n!$~\cite{KF95,STC06}.

A mapping $c: V(\Gamma) \rightarrow \{1,2,\ldots,k\}$ is called a {\it proper $k$--coloring} of a graph $\Gamma=(V,E)$ if $c(u)\neq c(v)$ whenever the vertices $u$ and $v$ are adjacent. The {\it chromatic number} $\chi(\Gamma)$ of a graph $\Gamma$ is the least number of colors needed to properly color vertices of $\Gamma$. A subset of vertices assigned to the same color forms an independent set, i.e. a proper $k$--coloring is the same as a partition of the vertex set into $k$ independent sets. The trivial lower and upper bounds on the chromatic number of the Pancake graphs are given as follows:
\begin{equation} \label{e1}
3\leqslant \chi(P_n)\leqslant n-1 {\mbox{ for any }}  n\geqslant 4.
\end{equation}

Indeed, the graph $P_n$ is $(n-1)$--regular, hence by Brooks' theorem~\cite{Br41} we have the upper bound. Moreover, $\chi(P_3)=2$ since $P_3\cong C_6$, and $\chi(P_4)=3$ since there are $7$--cycles in $P_n$ for any $n\geqslant 4$~\cite{KM10} which gives us the lower bound. The Brooks' bound is improved by $1$ for graphs with $\omega \leqslant (\Delta-1)/2$, where $\omega$ and $\Delta$ are the size of the maximum clique and the maximum degree of the graph (see~\cite{BK77,C78-1}). Since $\omega(P_n)=2$, then $\chi(P_n)\leqslant n-2$ for any $n\geqslant 6$. Moreover, there is a proper $3$--coloring of $P_5$~\cite{K17}. Thus, we have:
\begin{equation} \label{e2}
\chi(P_n)\leqslant n-2 {\mbox{ for any }}  n\geqslant 5.
\end{equation}

Catlin's bound for $C_4$--free graphs~\cite{C78-2}, that is $\chi\leqslant \frac{2}{3}\left(\Delta+3\right)$, gives one more bound for any $n\geqslant 8$:
\begin{equation} \label{e3}
\chi(P_n)\leqslant \frac{2}{3}\left(n + 2\right).
\end{equation}

Using structural properties of $P_n$, the following bounds were obtained in~\cite{K17}:

\begin{equation} \label{e4}
 {\mbox \it for } \ 5\leqslant n\leqslant 8, \ \chi(P_n)\leqslant \left\{
  \begin{array}{ll}
    n-k, & \hbox{if $n\equiv k\,(${\rm mod} $4)$ for $k=1,3$;} \\
    n-2, & \hbox{if $n$ is even;}
  \end{array}
\right.
\end{equation}

\begin{equation} \label{e5}
{\mbox \it for } \ 9\leqslant n\leqslant 16, \ \chi(P_n)\leqslant \left\{
  \begin{array}{ll}
    n-(k+2), & \hbox{if $n\equiv k\,(${\rm mod} $4)$ for $k=1,3$;} \\
    n-4, & \hbox{if $n$ is even;}
  \end{array}
\right.
\end{equation}

\begin{equation} \label{e6}
{\mbox \it for } \ n\geqslant 17, \ \chi(P_n)\leqslant \left\{
  \begin{array}{ll}
    n-(k+4), & \hbox{if $n\equiv k\,(${\rm mod} $4)$ for $k=1,2,3$;} \\
    n-8, & \hbox{if $n\equiv 0\,(${\rm mod} $4)$.}
  \end{array}
\right.
\end{equation}

These bounds improve~(\ref{e2}) for $n\geqslant 7$, however Catlin's bound~(\ref{e3}) is still better for all $n>28$ and some smaller $n$ (for example, $n=21,25,26,27$). Thus, they are far from good. Meanwhile, the asymptotic bound $\chi(P_n) \leqslant O\left(\frac{n-1}{log (n-1)}\right)$ holds for the Pancake graphs which follows from the results for $C_3, C_4$--free graphs~\cite{J96,K95}.

In this paper in Section~\ref{mr} we present a new upper bound which improves Catlin's bound~(\ref{e3}). The new bound is obtained using a subadditivity property of the chromatic number of $P_n$ and known chromatic numbers for $n\leqslant 9$. We have $\chi(P_3)=2$ since $P_3\cong C_6$, and $\chi(P_4)=3$ since there are $7$--cycles in $P_n,\,n\geqslant 4$. An example of a proper $3$--coloring for $P_5$ was given in~\cite{K17}. An optimal $4$-coloring for $P_6$ was computed by Toma\v{z} Pisanski, University of Primorska, Koper, Slovenia, and Jernej Azarija,  University of Ljubljana, Slovenia, so $\chi(P_6)=4$. Since $P_{n-1}$ is an induced subgraph of $P_n$, then $\chi(P_7)$ is at least $4$, and from~(\ref{e4}) we have $\chi(P_7)=4$. Optimal $4$-colorings for $P_8$ and $P_9$ were computed  by A.~H.~Ghodrati, Sharif University, Tehran, Iran. By~(\ref{e5}), $ 4 \leqslant \chi(P_n) \leqslant 12$, where $10 \leqslant n \leqslant 16$, however, proper $4$-colorings in these cases are unknown.  The known chromatic numbers are presented in the Table~\ref{chromN}.

In Section~\ref{cp} an equitable coloring is considered. A graph $\Gamma$ is said to be {\it equitably $k$-colorable} if $\Gamma$ has a proper $k$-coloring such that the sizes of any two color classes differ by at most one. The {\it equitable chromatic number} $\chi_{=}(\Gamma)$ is the smallest integer $k$ such that $\Gamma$ is equitably $k$-colorable. Equitable coloring was introduced by W.~Meyer in $1973$ due to scheduling problems~\cite{M73}. Moreover, it was conjectured that every connected graph with maximum degree $\Delta$ has an equitable coloring with $\Delta$ or fewer colors, with the exceptions of complete graphs and odd cycles. A strengthened version~\cite{CLW94} of the conjecture states that each such graph has an equitable coloring with exactly $\Delta$ colors, with one additional exception, a complete bipartite graph in which both sides of the bipartition have the same odd number of vertices. A survey on equitable colorings can be found in~\cite{L13}.

In Section~\ref{ec} an equitable $(n-1)$-coloring based on efficient dominating sets in the Pancake graphs $P_n, n>2$, is presented. Moreover, in Section~\ref{oc} simple optimal equitable $4$-colorings for $P_5, P_6$ and $P_7$ are described.

Let us note that any equitable coloring of $P_n$ with at most $n$ colors has the property that the sizes of all color classes are equal since every integer at most $n$ divides $n!$. Thus, we have a strongly equitable coloring~\cite{FKM17}.

Since equitable coloring is a proper coloring with an additional condition, the inequality $\chi(P_n) \leqslant \chi_{=}(P_n)$ holds for any $n\geqslant 2$. However, since all above optimal colorings are strongly equitable we have conjectured.

\begin{conjecture}
For any $n\geqslant 2$, $$\chi(P_n)=\chi_{=}(P_n).$$
\end{conjecture}

\begin{table}
\begin{center}
\begin{tabular}{|c|c|c|c|c|c|c|c|c|}
\hline $n$         & 2 & 3 & 4 & 5 & 6 & 7 & 8 & 9\\
\hline $\chi(P_n)$ & 2 & 2 & 3 & 3 & 4 & 4 & 4 & 4\\
\hline
\end{tabular}
\caption{Chromatic numbers of the Pancake graphs $P_n, \ 2 \leqslant n \leqslant 9$.}\label{chromN}
\end{center}
\end{table}

\section{Improved upper bound}\label{mr}

Our main result is given by the following theorem.

\begin{theorem} \label{MR} For any $n\geqslant 9,$ the following holds for the Pancake graph $P_n$:
\begin{equation} \label{e7}
\chi(P_n)\leqslant 4\left\lfloor\frac n9\right\rfloor+\chi\left(P_{n\pmod9}\right).
\end{equation}
\end{theorem}

To prove this result we need more notation. Let $[n]=\{1,2,\ldots,n\}$. We consider a permutation $\pi=[\pi_1 \pi_2 \ldots \pi_n]$ written as a string in one-line notation, where $\pi_i=\pi(i)$ for any $i\in [n]$. For $K\subset [n]$, let $P_{n,K}$ be the induced subgraph of $P_n$ whose vertex set consists of all permutations $\pi$ with $\pi_1\in K$. By the symmetry of $P_n$, for any $k$-element subset $K$ of $[n]$, the induced subgraph $P_{n,K}$ is isomorphic to $P_{n,[k]}$, which is abbreviated to $P_{n,k}$.

We define a map $f_{n,k}:P_{n,k}\to P_k$ by removing the elements that are not in $[k]$. For example, for $n=5$ and $k=3$, the vertex $[14352]$ of $P_{5,3}$ is mapped to the vertex $[132]$ of $P_3$. It is clear that this mapping is a graph homomorphism. Note that $f_{n,k}$ is surjective, but not necessarily an isomorphism. In fact, $P_{n,k}$ is not even connected unless $k=n$ or $n=2$.

Since an $r$-coloring of a graph $\Gamma$ is equivalent to a graph homomorphism from $\Gamma$ to the complete graph $K_r$, this property implies that
\begin{equation} \label{e8}
\chi(P_{n,k})\leqslant \chi(P_k).
\end{equation}
Since $P_{n,k}$ always contains a subgraph isomorphic to $P_k$ (e.~g. the subgraph of all vertices that end with $k+1,k+2,\ldots,n$) it even follows that $\chi(P_{n,k})=\chi(P_k)$.

One more useful property says that all fibers $f_{n+1,n}: P_{n+1,n}\to P_n$ are of size $n$, which means that $\left|f_{n+1,n}^{-1}(v)\right|=n$ for every $v \in V(P_n)$. Indeed, for any permutation of length $n$ one can insert $n+1$ at $n$ different positions: the first position is forbidden since the vertices of $P_{n+1,n}$ start with an element from $[n]$.

The following theorem immediately gives the general upper bound of~(\ref{e7}) with taking into account $\chi(P_9)=4$.

\begin{theorem}
The chromatic number of the Pancake graph is subadditive, i.~e.
\begin{equation} \label{e9}
\chi(P_{n+m})\leqslant\chi(P_n)+\chi(P_m)
\end{equation}
for all positive integers $n$ and $m$.
\end{theorem}

\proof
Let the vertices of $P_{n+m}$ be partitioned into sets $U$ and $W$ such that $U$ contains permutations whose first element is in $[n]$ and $W$ contains permutations whose first element is in $\{n+1,\ldots,n+m\}$. The subgraphs $P_{n+m,U}$ and $P_{n+m,W}$ are isomorphic to $P_{n+m,n}$ and $P_{n+m,m}$, respectively. Hence, by~(\ref{e8}) the graph $P_{n+m,U}$ is $\chi(P_n)$-colorable, and $P_{n+m,W}$ is $\chi(P_m)$-colorable. Using disjoint color sets on both subgraphs proves the desired inequality. \hfill $\square$

\section{Equitable coloring and optimal colorings}\label{cp}

It was shown in Introduction that there are optimal colorings of the Pancake graphs found by computer experiments. Such computations do not provide us any structural insight. In this section we consider colorings of the Pancake graphs $P_n, n\geqslant 3$, using their structural properties. More precisely, we present equitable colorings of $P_n$ in $(n-1)$ colors. 

An equitable coloring is not the same as a {\it perfect coloring} for which the multiset of colors of all neighbors of a vertex depends only on its own color~\cite{F07}. This type of coloring gives a partition known as an {\it equitable partition}~\cite{G97} which are used in algebraic combinatorics, graph theory and coding theory. In coding theory such kind of partitions are known as perfect codes~\cite{AAKh01,KMP02}. Some general information about equitable partitions can be found in~\cite{BCGG19}.

The notion of perfect codes was generalized to the Pancake graphs in a natural way in~\cite{DS03}. An independent set $D$ of vertices in a graph $\Gamma$ is an {\it efficient dominating set} (or $1$-perfect code) if each vertex not in $D$ is adjacent to exactly one vertex in $D$. There are $n$ efficient dominating sets in $P_n$~\cite{DS03,Q06} given by:
\begin{equation} \label{e10}
D_i=\{[i\,\pi_2\,\ldots\,\pi_n]\},
\end{equation}
where $\pi_k\in [n]\backslash\{i\},\,k\in [n]\backslash\{1\}, \, i\in [n]$.
It is obvious that $|D_{i_1}\bigcap D_{i_2}|=\emptyset$, $i_1, i_2\in [n], \ i_1\neq i_2$, which immediately gives a proper $n$-coloring. Moreover, this coloring is perfect and equitable.

To present a proper $(n-1)$-coloring of the Pancake graphs based on efficient dominating sets we need to define such sets for induced subgraphs $P_{n-1}$ of $P_n$.

Due to the hierarchical structure, for any $n\geqslant3$ the graph $P_n$ has $n$ copies of $P_{n-1}(i)$ with the vertex set $V_i=\{[\pi_1 \ldots \pi_{n-1} i]\}$, where $\pi_k\in [n]\backslash \{i\}, k\in [n-1]$, $|V_i|=(n-1)!$. Any two copies $P_{n-1}(i), P_{n-1}(j), i\neq j$, are connected by $(n-2)!$ edges $\{[i \pi_2 \ldots \pi_{n-1} j],\,[j \pi_{n-1} \ldots \pi_2 i]\}$, where $[i \pi_2 \ldots \pi_{n-1} j] r_n=[j \pi_{n-1} \ldots \pi_2 i]$. Prefix--reversals  $r_j, 2\leqslant j \leqslant n-1$, define internal edges in all $n$ copies $P_{n-1}(i)$, and the prefix--reversal $r_n$ defines external edges between copies.

Efficient dominating sets of $P_{n-1}(j), j\in [n]$, contain all permutations with the last element fixed to $j$ and the first element fixed~\cite{K13}, namely:

\begin{equation} \label{e11} D_i^j=\{[i\,\pi_2\,\ldots\,\pi_{n-1}\,j]\},
\end{equation}
where $i,j \in [n], i\neq j$, $\pi_k\in [n]\backslash\{i,\,j\},\,k \in [n] \backslash \{1,n\}$. 
For any $i \in [n]$, the sets~(\ref{e10}) and~(\ref{e11}) are given by the following obvious relationship:
\begin{equation} \label{e12}
D_i=\bigcup_{i=1, j\neq i}^n D_i^j.
\end{equation}

\subsection{Equitable $(n-1)$-coloring}\label{ec}

We now present an equitable $(n-1)$-coloring based on efficient dominating sets. 
Let
\begin{equation}\label{e13}
  D=\{D_i^j\st i,j\in[n],i\ne j\}, \ \ |D|=n\,(n-1),
\end{equation}
and
\begin{equation}\label{e14}
  D^j=\{D_i^j\st i\in[n],i\ne j\}, \ \ j\in[n], \ \ |D^j|=n-1,
\end{equation}
where
\begin{equation}\label{e15}
  |D_i^j|=(n-2)!.
\end{equation}

Note that $D$ partitions the vertices of $P_n$, and $D^j$ partitions the vertices of $P_n$ that end with $j$. We now define a graph $Q_n$ whose vertices are the elements of $D$, and $X,Y\in D$ are adjacent in $Q_n$ if and only if a vertex of $X$ is adjacent to a vertex of $Y$ in $P_n$. From the properties of the Pancake graphs we immediately see that vertices $D_i^j$ and $D_{i'}^{j'}$ are adjacent in $Q_n$ if and only if one of the following statements is true:
\begin{enumerate}
	\item [(A1)] $j=j'$ and $i\ne i'$.
	\item [(A2)] $i=j'$ and $j=i'$.
\end{enumerate}

It is obvious that a proper coloring $c$ of $Q_n$ trivially extends to a proper coloring of $P_n$ in such a way that any vertex of $P_n$ belongs to exactly one efficient dominating set $X\in D$ and we give it the color $c(X)$.

\begin{figure}[h]
     \centering
     \includegraphics[scale=0.65]{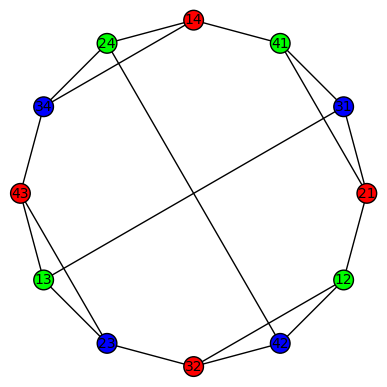}
     \caption{The equitable $3$-coloring of $Q_4$}\label{Q_4}
\end{figure}

\begin{figure}[h]
     \centering
     \includegraphics[scale=0.9]{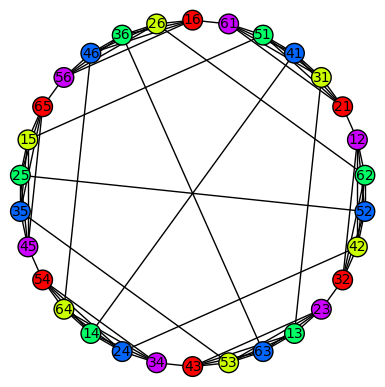}
     \caption{The equitable $5$-coloring of $Q_6$}\label{Q_6}
\end{figure}

We now have reduced the problem to finding a proper ($n-1$)-coloring for $Q_n$. First let us show an idea of such colorings for the graphs $Q_4$ and $Q_6$. The graphs are presented on Figures~\ref{Q_4} and~\ref{Q_6} such that the vertices corresponding to the set $D_i^j, \ i,j\in[n], i\ne j$, are denoted by labels $ij$. The vertices are arranged in a hamiltonian cycle such that all vertices with the same last element are grouped together and form $3$- and $5$-cliques, respectively. Within each clique the first element of labeling is cyclically incremented. Obviously, the elements of each clique must all have different colors, but the pictures suggest that we can {\lq}almost{\rq} cyclically repeat a color pattern chosen on the first clique. The only collisions occur with the {\lq}long{\rq} chords that connect antipodal vertices. We see that if we exchange the color of one end of each long chord with the color of the vertex counterclockwise next to it on the cycle, we obtain a proper coloring. It is clear that the proper colorings of $Q_4$ and $Q_6$ are equitable. Indeed, by~(\ref{e13})-(\ref{e15}) and from the construction $Q_4$ has $3$ color classes of cardinality $4$ each, and $Q_6$ has $5$ color classes of cardinality $6$ each. However, they are not perfect since the multisets of colors of all neighbors are different for different vertices having the same color. For example, in $Q_4$ the red vertex $(14)$ has two green and one blue neighbors, while the red vertex $(32)$ has two blue and one green neighbors. Similar, in $Q_6$ despite the red and the purple vertices have the same multiset of colors of their neighbors, the green, the blue and the dark blue vertices do not meet this condition to be perfect.

Note that this coloring is exactly the greedy coloring for the vertex sequence that starts with $(1n)$ and then counterclockwise follows the cycle.

\bs

Now we formalize and prove this observation. First we define a map:
\begin{equation} \label{e16}
f: D\to [n-1]
\end{equation}
such that
\begin{equation}\label{f-map}
\left\{
  \begin{array}{ll}
    f(D_i^j)=i-j, & \hbox{if $i>j$};\\
    f(D_i^j)=n+i-j, & \hbox{otherwise}.
  \end{array}
\right.
\end{equation}

Note that $f$ indeed has all its values in $[n-1]$, and that the restriction of $f$ to $D_j$ is injective, i.e.
$$f(D_i^j)=f(D_{i'}^j)\mbox{ if and only if }i=i'.$$

Next we let $k=\frac n2$ and define a coloring
\begin{equation} \label{e18}
c: D\to[n-1]
\end{equation}
by
\begin{equation}
c(D_i^j)=f(D_i^j)+\eps,
\end{equation}
where
\begin{equation}\label{eps}
\eps=
\left\{
  \begin{array}{ll}
    +1, & \hbox{if $j>k$ and $f(D_i^j)=k$;}\\
    -1, & \hbox{if $j>k$ and $f(D_i^j)=k+1$;}\\
     0, & \hbox{otherwise.}
  \end{array}
\right.
\end{equation}

If $n$ is odd, the first two cases cannot occur, since $f(D_i^j)$ is an integer, but $k$ is not. The restriction of $c$ to $D^j$ is still injective, since it is either equal to $f$ for $j\leqslant k$ or it is equal to $f$ with at most two adjacent function values exchanged.

In terms of the intuitive approach above, $f$ is the coloring that cyclically repeats along the cycle, and $c$ is the coloring that exchanges the colors near the ends of long chords.

\begin{theorem}
The coloring $c$ is proper and equitable for $P_n, n>2$.
\end{theorem}

\proof Indeed, suppose that $X=D_i^j$ and $Y=D_{i'}^{j'}$ are adjacent. If $j=j'$, then $i\ne i'$ by (A1). By the injectivity of $c$ for fixed $j$, $X$ and $Y$ must have different colors. By (A2) the only other possibility is that $Y=D_j^i$. Without loss of generality we assume that $j<i$.

First we handle the case that $n$ is odd. Then $c(X)=c(D_i^j)=f(D_i^j)=i-j$, and $c(Y)=c(D_j^i)=f(D_j^i)=n+j-i$. The equality $c(X)=c(Y)$ implies $n=2i-2j$, so $n$ is even. There is a contradiction.

In the last part of the proof we assume that $n$ is even, so $k$ is an integer. We do a case analysis based on the position of $k$.

If $j<i\leqslant k$ then $c(X)=i-j$, and $c(Y)=n+j-i$. The equality $c(X)=c(Y)$ implies $2i=n+2j=2k+2j$, so $i=k+j>k$, and this gives a contradiction.

If $j\leqslant k<i$ then $c(X)=i-j$, and $c(Y)=n+j-i+\eps$. The equality $c(X)=c(Y)$ implies $i-j=n+j-i+\eps$, so $2k+2j=2i-\eps$. Then $\eps$ must be even, hence $\eps=0$, and $i=k+j$ or $f(Y)=f(D_j^i)=n+j-i=k$. Since $i>k$ we have $\eps=1$ by~(\ref{eps}) which gives a contradiction.

If $k<j<i$ then $c(X)=i-j+\eps_1$ and $c(Y)=n+j-i+\eps_2$. The equality $c(X)=c(Y)$ implies $i-j+\eps_1=n+j-i+\eps_2$, so $2i+\eps_1=2k+2j+\eps_2$, and hence $\eps_1$ and $\eps_2$ have the same parity. There are the following possibilities.

(i) If $\eps_1=\eps_2$ then $i-j=k$. Since $k<j<i\leqslant n=2k$ we have $k=i-j<2k-k=k$ which gives a contradiction.

(ii) If $\eps_1=-1$, $\eps_2=1$. As in case~(i), we have $k+1=i-j<2k-k=k$ which leads to a contradiction.

(iii) If $\eps_1=1$, $\eps_2=-1$ then $i-j=k-1$. Again, since $k<j<i\leqslant 2k$, this is only possible if $j=k+1$ and $i=n=2k$. Then $f(D_i^j)=k-1$, so $\eps_1=0$ by~(\ref{eps}), and we have a contradiction.

In all cases the assumption of color equality leads to a contradiction, which finishes the proof that $c$ is a proper coloring. Moreover, by~(\ref{e18}) we have $(n-1)$ color classes each of which has a cardinality $|D^j|\cdot|D^j_i|=n(n-2)!$ (see~(\ref{e14}),~(\ref{e15})) which finishes the proof of the statement. \hfill $\square$

\subsection{Optimal colorings}\label{oc}

It is obvious that above equitable $(n-1)$-coloring produces an optimal coloring for $P_3$ and $P_4$ (see Table~\ref{chromN}). However, for $n>4$ a proper coloring of $Q_n$ can never produce an optimal coloring for $P_n$. Indeed, by~(\ref{e2}) for $n>4$ we have $\chi(P_n)<n-1$ and  $\chi(Q_n) \geqslant n-1$ since $Q_n$ contains an ($n-1$)-clique.

Now we give a simple optimal $4$-coloring of $P_5$, $P_6$ or $P_7$. We define an {\it even} ({\it odd}) prefix--reversal $r_i, \ 2\leqslant i \leqslant n$, if it corresponds to an even (odd) permutation. By~\cite[Lemma~4]{K17}, if $i\equiv0,1\pmod4$ then $r_i$ is an even prefix--reversal. Similar, $r_i$ is an odd prefix--reversal if $i\equiv2,3\pmod4$. By~\cite[Lemma~6]{K17}, the Pancake graph $P_n, n\geqslant 3$, has $n!/\ell$ independent even $\ell$--cycles where $6\leqslant \ell \leqslant 2n$.

Let $\Gamma$ be one of $P_5$, $P_6$, $P_7$. Then the subgraph $H$ generated by the even prefix--reversals $r_4$ and $r_5$ is a spanning subgraph of $\Gamma$ consisting of disjoint $10$-cycles $C_{10}=(r_5\ r_4)^5$. Since even prefix--reversals preserve parity, all vertices of $H$ that are on the same cycle have the same parity. Since all other edges of $\Gamma$ correspond to odd permutations, they have one endpoint on an {\lq}even{\rq} $10$-cycle and the other endpoint on an {\lq}odd{\rq} $10$-cycle. Thus, we can $2$-color the even cycles and $2$-color the odd cycles using two other colors. This results in an equitable $4$-coloring of $\Gamma$ where each of the color classes has $n!/4$ vertices. Since $\chi(P_6)=\chi(P_7)=4$ this coloring is optimal for $P_6$ and $P_7$.

\section{Discussion}

There are trivial examples of graphs for which Conjecture $1$ holds such as even cycles, bipartite graphs with equal parts. Any Cayley graph over the symmetric group generated by a set of transpositions gives a bipartite graph with two equal color classes. The statement holds for the Hamming graphs $H(d,q)$ whose (equitable) chromatic number is $q$. A regular graph with a Hoffman coloring always gives a strongly equitable coloring~\cite{BBH07}. Hoffman's lower bound is known as $\chi(\Gamma)\geqslant 1-\lambda_1/\lambda_v$, where $\lambda_1$ and $\lambda_v$ are the largest and the smallest eigenvalues of $\Gamma$. If equality holds, an optimal coloring of $\Gamma$ is called a Hoffman coloring.

\section*{Acknowledgements} The research work of the second author is supported by Mathematical Center in Akademgorodok, the agreement with Ministry of Science and High Education of the Russian Federation number 075-15-2019-1613. The authors thank Alexander Kostochka and Sergey Avgustinovich for useful discussions on equitable colorings.

\ms
\end{document}